\documentclass[12pt,a4paper,fleqn]{article}
\usepackage{a4wide,amsfonts,amsmath,latexsym,amssymb,euscript,graphicx,units,mathrsfs}

\usepackage{graphicx}
\usepackage{color}
\usepackage{amssymb}
\usepackage{amssymb}
\usepackage[T1]{fontenc}
\usepackage{latexsym}
\usepackage{xypic}
\usepackage{eufrak}
\usepackage{euscript}
\usepackage{amsfonts,amsmath}
\usepackage{verbatim}
\usepackage[english]{babel}
\usepackage{mathrsfs}
\usepackage{units}

\newtheorem{prop}{Proposition}[section]

\newtheorem{lemma}[prop]{Lemma}

\newtheorem{remark}[prop]{Remark}
\newtheorem{thm}[prop]{Theorem}

\newtheorem{defi}[prop]{Definition}

\renewcommand{\geq}{\geqslant}
\def\leq{\leqslant}

\newcommand{\N}{\mathbb{N}}

\newcommand{\R}{\mathbb{R}}

\def\HH{\EuFrak H}

\def\1{{\mathbf{1}}}

\def\1{{\mathbf{1}}}
\def\0.5{{\frac{1}{2}}}

\newcommand{\fin}
{ \vspace{-0.6cm}
\begin{flushright}
\mbox{$\Box$}
\end{flushright}
\noindent }

\newcommand{\qed}{\nopagebreak\hspace*{\fill}
{\vrule width6pt height6ptdepth0pt}\par}

\begin{document}

\begin{center}
{\large{\bf On the Gaussian approximation
of vector-valued multiple integrals}}\\~\\
Salim Noreddine\footnote{
Laboratoire de Probabilit{\'e}s et
Mod{\`e}les Al{\'e}atoires,
Universit{\'e} Paris 6,
Bo{\^\i}te courrier 188, 4 Place Jussieu, 
75252 Paris Cedex 5, France.
Email: {\tt salim.noreddine@polytechnique.org}
}
and
Ivan Nourdin\footnote{Institut \'Elie Cartan, Universit\'e Henri Poincar\'e, BP 239, 54506 Vandoeuvre-l\`es-Nancy, France. Email: {\tt inourdin@gmail.com}}\footnote{Supported in part by the (french) ANR grant `Exploration des Chemins Rugueux'}\\
{\it Université Paris 6} and {\it Universit\'e Nancy 1}\\~\\
\end{center}
{\small \noindent {\bf Abstract:} 
By combining the findings of two recent, seminal papers by Nualart, Peccati and Tudor, we get that 
the convergence in law of any sequence of vector-valued multiple integrals $F_n$ towards
a centered Gaussian random vector $N$, 
with given covariance matrix $C$,
is reduced to just the convergence of: $(i)$ the fourth
cumulant of each component of $F_n$ to zero; $(ii)$
the covariance matrix of $F_n$ to $C$.
The aim of this paper is to understand more deeply this somewhat surprising phenomenom.
To reach this goal, we offer two results of different nature. The first one is
an explicit bound for $d(F,N)$ in terms of the fourth cumulants of the components of $F$, when
$F$ is a $\R^d$-valued 
random vector whose components are multiple integrals of possibly different orders, $N$ is
the Gaussian counterpart of $F$ (that is, a Gaussian centered vector
sharing the same covariance with $F$)
and $d$ stands for the Wasserstein distance.
The second one is a new expression for the cumulants of $F$ as above, from which
it is easy to derive yet another proof of the previously quoted result by Nualart, Peccati and Tudor.\\

\noindent {\bf Keywords:} Central limit theorem; Cumulants; Malliavin calculus; Multiple integrals; 
Ornstein-Uhlenbeck Semigroup.\\

\noindent
{\bf 2000 Mathematics Subject Classification:} 60F05; 60G15; 60H05; 60H07. }
\\

\section{Introduction}\label{intro}
Let $B=(B_t)_{t\in[0,T]}$ be a standard Brownian motion.
The following result, proved in \cite{NP,PT}, yields a very surprising condition
under which a sequence of vector-valued multiple integrals converges in law to a Gaussian random vector.
(If needed, we refer the reader to section \ref{S : Malliavin} for the exact meaning of 
$\int_{[0,T]^q}f(t_1,\ldots,t_q)dB_{t_1}\ldots dB_{t_q}$.)
\begin{thm}[Nualart-Peccati-Tudor]\label{npt-thm}
Let  
$q_d,\ldots,q_1\geq 1$ be some fixed integers. Consider a $\R^d$-valued random sequence
of the form
\begin{eqnarray*}
F_n&=&(F_{1,n},\ldots,F_{d,n})\\
&=&\left(
\int_{[0,T]^{q_1}}f_{1,n}(t_1,\ldots,t_{q_1})dB_{t_1}\ldots dB_{t_{q_1}},
\ldots,
\int_{[0,T]^{q_d}}f_{d,n}(t_1,\ldots,t_{q_d})dB_{t_1}\ldots dB_{t_{q_d}}
\right),
\end{eqnarray*}
where each $f_{i,n}\in L^2([0,T]^{q_i})$ is supposed to be symmetric.
Let $N\sim\mathscr{N}_d(0,C)$ be a centered Gaussian random vector on $\R^d$ with covariance matrix $C$.
Assume furthermore that 
\begin{equation}\label{cij}
\lim_{n\to\infty}E[F_{i,n}F_{j,n}]= C_{ij}\quad\mbox{for all $i,j=1,\ldots,d$.} 
\end{equation}
Then, as $n\to\infty$, the following two assertions are equivalent:
\begin{enumerate}
\item[(i)] $F_n\overset{{\rm Law}}{\longrightarrow} N$;
\item[(ii)] $\forall i=1,\ldots,d$: \,$E[F_{i,n}^4]- 3E[F_{i,n}^2]^2\to 0$.
\end{enumerate}
\end{thm}

This theorem
represents a drastic simplification with respect to the method of moments. 
The original proofs performed in \cite{NP,PT} are both based on 
tools coming from Brownian stochastic analysis, such as the Dambis, Dubins and Schwarz theorem.
In \cite{NuOL}, Nualart and Ortiz-Latorre gave an alternative proof
exclusively using the basic
operators $\delta$, $D$ and $L$ of Malliavin calculus.
Later on, combining Malliavin calculus with  
Stein's method in the spirit
of \cite{stein},
Nourdin, Peccati and R\'eveillac were able to
associate 
an explicit bound to convergence $(i)$ in Theorem \ref{npt-thm}:
\begin{thm}[see \cite{NPReveillac}]\label{npr-thm}
Consider a $\R^d$-valued random vector
of the form
\begin{eqnarray*}
F&=&(F_{1},\ldots,F_{d})\\
&=&\left(
\int_{[0,T]^{q_1}}f_{1}(t_1,\ldots,t_{q_1})dB_{t_1}\ldots dB_{t_{q_1}},
\ldots,
\int_{[0,T]^{q_d}}f_{d}(t_1,\ldots,t_{q_d})dB_{t_1}\ldots dB_{t_{q_d}}
\right),
\end{eqnarray*}
where $q_1,\ldots,q_d\geq 1$ are some given integers, and
each $f_{i}\in L^2([0,T]^{q_i})$ is symmetric.
Let $C=(C_{ij})_{1\leq i,j\leq d}\in\mathcal{M}_d(\R)$ be the covariance matrix
of $F$, i.e. $C_{ij}=E[F_iF_j]$. Consider
a centered Gaussian random vector 
$N\sim\mathscr{N}_d(0,C)$ with same covariance matrix $C$.
Then:
\begin{equation}\label{def-wass}
d_{1}(F,N):=\sup_{h\in {\rm Lip}(1)}\big| E[h(F)] - E[h(N)]\big|
\leq \| C^{-1}\|_{op} \, \| C\|_{op}^{1/2}\,\Delta_C(F),
\end{equation}
with the convention $\|C^{-1}\|_{op}=+\infty$ whenever $C$ is not invertible.
Here:

- ${\rm Lip}(1)$ is the set of Lipschitz functions with constant 1 (that is,
the set of functions $h:\R^d\to\R$ so that $|h(x)-h(y)|\leq \|x-y\|_{\R^d}$ for all $x,y\in\R^d$),\\

- $\|C\|_{op}=\sup_{x\in\R^d\setminus\{0\}}\|Cx\|_{\R^d}/\|x\|_{\R^d}$ denotes the operator norm 
on $\mathcal{M}_d(\R)$,\\

- the quantity $\Delta_C(F)$ is defined as
\begin{equation}\label{delta}
\Delta_C(F):=\sqrt{\sum_{i,j=1}^d
E\left[\left( C_{ij} - \frac1{q_j}\langle DF_{i}, DF_{j} \rangle_{L^2([0,T])}\right)^2\right] },
\end{equation}
where $D$ indicates the Malliavin derivative operator (see section \ref{S : Malliavin})
and $\langle\cdot,\cdot\rangle_{L^2([0,T])}$ is the usual inner product on $L^2([0,T])$.
\end{thm}
When the covariance matrix $C$ of $F$
is not invertible (or when one is not able to check whether it is or not), one is forced to work 
with functions $h$ that are smoother than the one involved in the definition
(\ref{def-wass})
of $d_1(F,N)$. To this end, we adopt the following simplified notation
for functions $h:\R^d\to\R$ belonging to $\mathcal{C}^2$:
\begin{equation}\label{h''}
\|h''\|_\infty = \max_{i,j=1,\ldots,d}\,\,\sup_{x\in\R^d}\left|\frac{\partial^2h}{\partial x_i\partial x_j}(x)\right|.
\end{equation}
\begin{thm}[see \cite{NP-survey}]\label{np-thm}
Let the notation and assumptions of Theorem \ref{npr-thm} prevail. Then:
\begin{equation}\label{def-smooth}
d_{2}(F,N):=\sup_{\|h''\|_\infty\leq 1}\big| E[h(F)] - E[h(N)]\big|
\leq \frac12 \,\Delta_C(F),
\end{equation}
with $\Delta_C(F)$ still given by (\ref{delta}).
\end{thm}

Are the upper bounds (\ref{def-wass})-(\ref{def-smooth}) in Theorems \ref{npr-thm} and \ref{np-thm} relevant? 
Yes, very! 
Indeed, we have the following proposition.
\begin{prop}[see \cite{NuOL}]\label{nol-lem}
Let the notation and assumptions of Theorem \ref{npt-thm} prevail. Recall the definition
(\ref{delta}).
Then, as $n\to\infty$,
$\Delta_C(F_n)\to 0$ if and only if $E[F_{i,n}^4] - 3E[F_{i,n}^2]^2\to 0$ for all $i$.
\end{prop}
In the present paper, as a first result we offer the following quantitative
version of Proposition \ref{nol-lem}.

\begin{thm}\label{main1}
Let the notation and assumptions of Theorem \ref{npr-thm} prevail, and recall
the definition (\ref{delta}) of $\Delta_C(F)$.
Then:
\begin{equation}\label{1.7}
\Delta_C(F)\leq 
\psi\!\left(E[F_1^4]-3E[F_1^2]^2,E[F_1^2],\ldots,E[F_d^4]-3E[F_d^2]^2,E[F_d^2]
\right),
\end{equation}
with $\psi:(\R\times\R_+)^d\to\R$ the function defined as
\begin{eqnarray}
\psi\big(x_1,y_1,\ldots,x_d,y_d\big)\!\!&=&\!\!\!
\sum_{i,j=1}^d 
{\bf 1}_{\{q_i=q_j\}}\sqrt{2\sum_{r=1}^{q_i-1}\binom{2r}{r}}|x_i|^{1/2}
+\sum_{i,j=1}^d {\bf 1}_{\{q_i\neq q_j\}}\bigg\{\sqrt2
\sqrt{y_j}|x_i|^{1/4}\notag\\
&&+
\left.\sum_{r=1}^{q_i\wedge q_j-1} \sqrt{2(q_i+q_j-2r)!}\binom{q_j}{r}|x_i|^{1/2}
\right\}.
\label{quantitativeNP}
\end{eqnarray}
\end{thm}
Since, for each compact $B\subset(0,\infty)^d$, it is readily checked that there exists a constant $c_{B,q_1,\ldots,q_d}>0$ so that
\[
\sup_{(y_1,\ldots,y_d)}\psi(x_1,y_1, \ldots,x_s,y_d)\leq
c_{B,q_1,\ldots,q_d}\,\sum_{i=1}^d \left(|x_i|^{1/4}+|x_i|^{1/2}\right),
\]
we immediately see that
the upper bound 
(\ref{1.7}), together with Theorem \ref{np-thm}, now show in a
clear manner why $(ii)$ implies $(i)$ in Theorem \ref{npt-thm}. \\

In a second part of this paper,
we are interested in `calculating', by means of the basic operators $D$ and $L$
of Malliavin calculus, the cumulants of any
vector-valued functional $F$ of the Brownian motion $B$.
(Actually, we will even do so for functionals of any given {\it isonormal Gaussian process} $X$).
In fact, this part is nothing but the multivariate extension of the results obtained by Nourdin and Peccati
in \cite{NP-cumulants}.

Then, in the particular case where the components of $F$ 
have the form of a multiple Wiener-It\^{o} integral (as in
Theorem \ref{npr-thm}), 
our formula leads to a new compact representation for 
the cumulants of $F$ (see Theorem \ref{thm1.6} just below), implying in turn 
yet another proof of Theorem \ref{npt-thm} (see section \ref{sec4.3}).
 
\begin{thm}\label{thm1.6}
Let $m\in\N^d\setminus\{0\}$ with $|m|\geq 3$. Write $m=l_1+\ldots+l_{|m|}$,
where $l_i\in\{e_1,\ldots,e_d\}$ for each $i$.
(Up to possible permutations of factors, we have existence and uniqueness
of this decomposition of $m$.)
Consider a $\R^d$-valued random vector
of the form
\begin{eqnarray*}
F&=&(F_{1},\ldots,F_{d})\\
&=&\left(
\int_{[0,T]^{q_1}}f_{1}(t_1,\ldots,t_{q_1})dB_{t_1}\ldots dB_{t_{q_1}},
\ldots,
\int_{[0,T]^{q_d}}f_{d}(t_1,\ldots,t_{q_d})dB_{t_1}\ldots dB_{t_{q_d}}
\right),
\end{eqnarray*}
where $q_1,\ldots,q_d\geq 1$ are some given integers, and
each $f_{i}\in L^2([0,T]^{q_i})$ is symmetric.
When $l_k=e_j$, we set $\lambda_k=j$, so that
$F^{l_k}=F_{\lambda_k}$ for all $k=1,\ldots,|m|$.
Then:
\[
\kappa_m(F)= (q_{\lambda_{|m|}})!(|m|-1)!\sum c_{q,l}(r_2,\ldots,r_{|m|-1})
\langle
f_{\lambda_1}\widetilde{\otimes}_{r_2} f_{\lambda_2}\ldots \widetilde{\otimes}_{r_{{|m|-1}}}f_{\lambda_{|m|-1}},
f_{\lambda_{|m|}}
\rangle_{
L^2([0,T]^{q_{\lambda_{|m|}}})
},
\]
where the sum $\sum$ runs over all collections of integers $r_2,\ldots,r_{|m|-1}$ such that:
\begin{enumerate}
\item[(i)] $1\leq r_i\leq q_{\lambda_i}$ for all $i=2,\ldots,|m|-1$;
\item[(ii)] $r_2+\ldots+r_{|m|-1}=\frac{
q_{\lambda_1}+\ldots+q_{\lambda_{|m|-1}}-q_{\lambda_{|m|}}}{2}$;
\item[(iii)] $r_2 <\frac{q_{\lambda_1}+q_{\lambda_2}}{2}$, $\ldots$, 
$r_2+\ldots +r_{|m|-2}< \frac{q_{\lambda_1}+\ldots+q_{\lambda_{|m|-2}}}2$;
\item[(iv)] $r_3\leq q_{\lambda_1}+q_{\lambda_2}-2r_2$, 
$\ldots$, 
$r_{|m|-1}\leq q_{\lambda_1}+q_{\lambda_{|m|-2}}-2r_2-\ldots-2r_{|m|-2}$;
\end{enumerate}
and where the combinatorial constants $c_{q,l}(r_2,\ldots,r_{s})$ are recursively defined by the
relations
\[
c_{q,l}(r_2)=q_{\lambda_2}(r_2-1)!\binom{q_{\lambda_1}-1}{r_2-1}\binom{q_{\lambda_2}-1}{r_2-1},
\]
and, for $s\geq 3$,
\begin{eqnarray*}
c_{q,l}(r_2,\ldots,r_{s})&=&q_{\lambda_s}(r_{s}-1)!\binom{q_{\lambda_1}+\ldots+q_{\lambda_s}-2r_2-\ldots - 
2r_{s-1}-1}{r_{s}-1}\\
&&\hskip6.5cm\times
\binom{q_{\lambda_s}-1}{r_{s}-1}c_{q,l}(r_2,\ldots,r_{s-1}).
\end{eqnarray*}
\end{thm}

The rest of the paper is organized as follows.
Section 2 gives (concise) background and notation for Malliavin calculus.
The proof of Theorem \ref{main1} is performed in Section 3.
Finally, Section 4 is devoted to the study of cumulants,
and contains in particular the proof of Theorem \ref{thm1.6}.

\section{Preliminaries on Malliavin calculus}\label{S : Malliavin}
In this section, we present the basic elements of Gaussian analysis and Malliavin calculus 
that are used throughout this paper. 
The reader
is referred to \cite{Nbook} for any unexplained definition or
 result.

\smallskip

Let $\EuFrak H$ be a real separable Hilbert space. For any $q\geq 1$, 
let $\EuFrak H^{\otimes q}$ be the $q$th tensor power of $\EuFrak H$, and denote
by $\EuFrak H^{\odot q}$ the associated $q$th symmetric tensor power. We write $X=\{X(h),h\in \EuFrak H\}$ to indicate
an isonormal Gaussian process over
$\EuFrak H$ (fixed once for all), defined on some probability space $(\Omega ,\mathcal{F},P)$.
This means that $X$ is a centered Gaussian family, whose covariance is given by the relation
$E\left[ X(h)X(g)\right] =\langle h,g\rangle _{\EuFrak H}$. We also assume that $\mathcal{F}=\sigma(X)$, that is,
$\mathcal{F}$ is generated by $X$.

For every $q\geq 1$, let $\mathcal{H}_{q}$ be the $q$th Wiener chaos of $X$,
defined as the closed linear subspace of $L^2(\Omega ,\mathcal{F},P)$
generated by the family $\{H_{q}(X(h)),h\in \EuFrak H,\left\|
h\right\| _{\EuFrak H}=1\}$, where $H_{q}$ is the $q$th Hermite polynomial
given by \[
H_q(x) = (-1)^q e^{\frac{x^2}{2}}
 \frac{d^q}{dx^q} \big( e^{-\frac{x^2}{2}} \big).
\]
We write by convention $\mathcal{H}_{0} = \mathbb{R}$. For
any $q\geq 1$, the mapping $I_{q}(h^{\otimes q})=q!H_{q}(X(h))$ can be extended to a
linear isometry between the symmetric tensor product $\EuFrak H^{\odot q}$
(equipped with the modified norm $\sqrt{q!}\left\| \cdot \right\| _{\EuFrak H^{\otimes q}}$)
and the $q$th Wiener chaos $\mathcal{H}_{q}$. For $q=0$, we write $I_{0}(c)=c$, $c\in\mathbb{R}$. 
For $q=1$, we have $I_1(h)=X(h)$, $h\in\HH$. Moreover,
a random variable of the type $I_q(h)$, $h\in\HH^{\odot q}$,
has finite moments of all orders.
 
In the particular case where $\HH=L^2([0,T])$, one has that $(B_t)_{t\in[0,T]}=\big(X({\bf 1}_{[0,t]})\big)_{t\in[0,T]}$ 
is a standard Brownian motion. Moreover, $\HH^{\odot q}=L^2_s([0,T]^q)$ is
the space of symmetric and square integrable functions on $[0,T]^q$, and
\[
I_q(f)=:\int_{[0,T]^q}f(t_1,\ldots,t_q)dB_{t_1}\ldots dB_{t_q},\quad f\in\HH^{\odot q},
\] 
coincides with the multiple Wiener-It\^o integral of order $q$
of $f$ with respect to $B$, see \cite{Nbook} for further details about this point.

It is
well-known that $L^2(\Omega):=L^2(\Omega ,\mathcal{F},P)$
can be decomposed into the infinite orthogonal sum of the spaces $\mathcal{H}_{q}$. It follows that any square integrable random variable
$F\in L^2(\Omega)$ admits the following so-called chaotic expansion:
\begin{equation}
F=\sum_{q=0}^{\infty }I_{q}(f_{q}),  \label{E}
\end{equation}
where $f_{0}=E[F]$, and the $f_{q}\in \EuFrak H^{\odot q}$, $q\geq 1$, are
uniquely determined by $F$. For every $q\geq 0$, we denote by $J_{q}$ the
orthogonal projection operator on the $q$th Wiener chaos. In particular, if
$F\in L^2(\Omega)$ is as in (\ref{E}), then
$J_{q}F=I_{q}(f_{q})$ for every $q\geq 0$.

Let $\{e_{k}\}_{k\geq 1}$ be a complete orthonormal system in $\EuFrak H$.
Given $f\in \EuFrak H^{\odot p}$ and $g\in \EuFrak H^{\odot q}$, for every
$r=0,\ldots ,p\wedge q$, the contraction of $f$ and $g$ of order $r$
is the element of $\EuFrak H^{\otimes (p+q-2r)}$ defined by
\begin{equation}
f\otimes _{r}g=\sum_{i_{1},\ldots ,i_{r}=1}^{\infty }\langle
f,e_{i_{1}}\otimes \ldots \otimes e_{i_{r}}\rangle _{\EuFrak H^{\otimes
r}}\otimes \langle g,e_{i_{1}}\otimes \ldots \otimes e_{i_{r}}
\rangle_{\EuFrak H^{\otimes r}}.  \label{v2}
\end{equation}
Note that the definition of $f\otimes_r g$ does not depend
on the particular choice of $\{e_k\}_{k\geq 1}$, and that
$f\otimes _{r}g$ is not necessarily symmetric; we denote its
symmetrization by $f\widetilde{\otimes }_{r}g\in \EuFrak H^{\odot (p+q-2r)}$.
Moreover, $f\otimes _{0}g=f\otimes g$ equals the tensor product of $f$ and
$g$, whereas $f\otimes _{q}g=\langle f,g\rangle _{\EuFrak H^{\otimes q}}$ whenever $p=q$.

It can be shown that the following product formula holds: if $f\in \EuFrak
H^{\odot p}$ and $g\in \EuFrak
H^{\odot q}$, then
\begin{eqnarray}\label{multiplication}
I_p(f) I_q(g) = \sum_{r=0}^{p \wedge q} r! {p \choose r}{ q \choose r} I_{p+q-2r} (f\widetilde{\otimes}_{r}g).
\end{eqnarray}
\smallskip

We now introduce some basic elements of the Malliavin calculus with respect
to the isonormal Gaussian process $X$. Let $\mathcal{S}$
be the set of all
cylindrical random variables of
the form
\begin{equation}
F=g\left( X(\phi _{1}),\ldots ,X(\phi _{n})\right) ,  \label{v3}
\end{equation}
where $n\geq 1$, $g:\mathbb{R}^{n}\rightarrow \mathbb{R}$ is an infinitely
differentiable function such that its partial derivatives have polynomial growth, and each
$\phi _{i}$ belongs to $\EuFrak H$.
The Malliavin derivative of $F$ with respect to $X$ is the element of
$L^2(\Omega ,\EuFrak H)$ defined as
\begin{equation*}
DF\;=\;\sum_{i=1}^{n}\frac{\partial g}{\partial x_{i}}\left( X(\phi
_{1}),\ldots ,X(\phi _{n})\right) \phi _{i}.
\end{equation*}
In particular, $DX(h)=h$ for every $h\in \EuFrak H$. By iteration, one can
define the $m$th derivative $D^{m}F$, which is an element of $L^2(\Omega ,\EuFrak H^{\odot m})$,
for every $m\geq 2$.
For $m\geq 1$ and $p\geq 1$, ${\mathbb{D}}^{m,p}$ denotes the closure of
$\mathcal{S}$ with respect to the norm $\Vert \cdot \Vert _{m,p}$, defined by
the relation
\begin{equation*}
\Vert F\Vert _{m,p}^{p}\;=\;E\left[ |F|^{p}\right] +\sum_{i=1}^{m}E\left[
\Vert D^{i}F\Vert _{\EuFrak H^{\otimes i}}^{p}\right] .
\end{equation*}
One also writes $\mathbb{D}^{\infty} = \bigcap_{m\geq 1}
\bigcap_{p\geq 1}\mathbb{D}^{m,p}$. The Malliavin derivative $D$ obeys the following chain rule. If
$\varphi :\mathbb{R}^{n}\rightarrow \mathbb{R}$ is continuously
differentiable with bounded partial derivatives and if $F=(F_{1},\ldots
,F_{n})$ is a vector of elements of ${\mathbb{D}}^{1,2}$, then $\varphi
(F)\in {\mathbb{D}}^{1,2}$ and
\begin{equation}\label{chainrule}
D\,\varphi (F)=\sum_{i=1}^{n}\frac{\partial \varphi }{\partial x_{i}}(F)DF_{i}.
\end{equation}
The conditions imposed on $\varphi$ for (\ref{chainrule}) to hold 
(that is, the partial derivatives of $\varphi$ must be bounded)
are by no means optimal. For instance, the chain rule combined with a classical
approximation argument leads to
$D(X(h)^m)=mX(h)^{m-1}h$ for $m\geq 1$ and $h\in\HH$.

We denote by $\delta $ the adjoint of the operator $D$, also called the
divergence operator. A random element $u\in L^2(\Omega ,\EuFrak H)$
belongs to the domain of $\delta $, noted $\mathrm{Dom}\delta $, if and
only if it verifies
$|E\langle DF,u\rangle _{\EuFrak H}|\leq c_{u}\,\Vert F\Vert _{L^2(\Omega)}$
for any $F\in \mathbb{D}^{1,2}$, where $c_{u}$ is a constant depending only
on $u$. If $u\in \mathrm{Dom}\delta $, then the random variable $\delta (u)$
is defined by the duality relationship
\begin{equation}
E[F\delta (u)]=E\langle DF,u\rangle _{\EuFrak H},  \label{ipp}
\end{equation}
which holds for every $F\in {\mathbb{D}}^{1,2}$.

The operator $L$ is defined as
$L=\sum_{q=0}^{\infty }-qJ_{q}$.
The domain of $L$ is
\begin{equation*}
\mathrm{Dom}L=\{F\in L^2(\Omega ):\sum_{q=1}^{\infty }q^{2}E[(J_qF)^2]<\infty \}=\mathbb{D}^{2,2}\text{.}
\end{equation*}
There is an important relation between the operators $D$, $\delta $ and $L$.
A random variable $F$ belongs to
$\mathbb{D}^{2,2}$ if and only if $F\in \mathrm{Dom}\left( \delta D\right) $
(i.e. $F\in {\mathbb{D}}^{1,2}$ and $DF\in \mathrm{Dom}\delta $) and, in
this case,
\begin{equation}
\delta DF=-LF.  \label{k1}
\end{equation}

For any $F \in L^2(\Omega )$, we define $L^{-1}F =\sum_{q=1}^{\infty }-\frac{1}{q} J_{q}(F)$. The operator $L^{-1}$ is called the
pseudo-inverse of $L$. Indeed, for any $F \in L^2(\Omega )$, we have that $L^{-1} F \in  \mathrm{Dom}L
= \mathbb{D}^{2,2}$,
and
\begin{equation}\label{Lmoins1}
LL^{-1} F = F - E[F].
\end{equation}
We end up these preliminaries on Malliavin calculus by stating a useful lemma, that
is going to be intensively used in 
the forthcoming Section \ref{multivariate}.
\begin{lemma}
Suppose that $F\in\mathbb{D}^{1,2}$ and $G\in L^2(\Omega)$. Then, $L^{-1}G\in\mathbb{D}^{2,2}$ and we have:
\begin{equation}\label{NPlemma}
E[FG]=E[F]E[G]+E[\langle DF,-DL^{-1}G\rangle_\HH].
\end{equation}
\end{lemma}
{\it Proof}. By (\ref{k1}) and (\ref{Lmoins1}), 
\[
E[FG]-E[F]E[G] = E[F(G-E[G])] = E[F\times LL^{-1}G] = E[F\delta(-DL^{-1}G)],
\]
and the result is obtained by using the integration by parts formula (\ref{ipp}).\qed

\section{Proof of Theorem \ref{main1}}\label{proof-quantitativeNP}

The aim of this section is to prove Theorem \ref{main1}.
We restate it here for convenience, by reformulating it in the more general context
of isonormal Gaussian process rather than Brownian motion.\\
\\
{\bf Theorem \ref{main1}}
{\it 
Let $X=\{X(h),\,h\in\HH\}$ be an isonormal Gaussian process, and 
$q_d,\ldots,q_1\geq 1$ be some fixed integers. 
Consider a $\R^d$-valued random vector
of the form
\begin{eqnarray*}
F=(F_{1},\ldots,F_{d})=\big(I_{q_1}(f_{1}),\ldots,I_{q_d}(f_{d})\big),
\end{eqnarray*}
where each $f_{i}$ belongs to $\HH^{\odot{q_i}}$.
Let $C=(C_{ij})_{1\leq i,j\leq d}\in\mathcal{M}_d(\R)$ be the covariance matrix
of $F$, i.e. $C_{ij}=E[F_iF_j]$, and consider
a centered {\rm Gaussian} random vector 
$N\sim\mathscr{N}_d(0,C)$ with same covariance matrix $C$.
Then
\begin{equation}\label{3.7}
\Delta_C(F)\leq 
\psi\!\left(E[F_1^4]-3E[F_1^2]^2,E[F_1^2],\ldots,E[F_d^4]-3E[F_d^2]^2,E[F_d^2]
\right),
\end{equation}
with $\Delta_C(F)$ given by (\ref{delta}), and
where $\psi:(\R\times\R_+)^d\to\R$ is the function given by (\ref{quantitativeNP}).\\
}

In order to prove Theorem \ref{main1}, we first need to gather
several results from the existing literature.
We collect them in the following lemma. We freely use the definitions and notation
introduced in sections \ref{intro} and \ref{S : Malliavin}.

\begin{lemma}
Let $F=I_p(f)$ and $G=I_q(g)$, with $f\in\HH^{\odot p}$ and
$g\in\HH^{\odot q}$ ($p,q\geq 1$).

1. If
$p=q$, one has the estimate:
\begin{eqnarray}
&&E\left[\left(E[FG]-\frac1p\left\langle DF,DG\right\rangle_\HH\right)^2\right]\label{ineq1}
\\
&&\hskip1cm\leq
\frac{p^2}{2}\sum_{r=1}^{p-1}(r-1)!^2\binom{p-1}{r-1}^4(2p-2r)!
\big( \|f\otimes_{p-r}f\|^2_{\HH^{\otimes
2r}}+\|g\otimes_{p-r}g\|^2_{\HH^{\otimes 2r}}\big),\notag
\end{eqnarray}
whereas, if $p< q$, one has that
\begin{eqnarray}
&&E\left[\left(\frac1q\left\langle DF,DG\right\rangle_\HH
\right)^2\right]\label{ineq2}
\leq p!^2\binom{q-1}{p-1}^2(q-p)!\|f\|^2_{\HH^{\otimes p}}\|g\otimes_{q-p}g\|_{\HH^{\otimes 2p}}\\
&&+
\frac{p^2}{2}\sum_{r=1}^{p-1}(r-1)!^2\binom{p-1}{r-1}^2\binom{q-1}{r-1}^2(p+q-2r)!\big(
\|f\otimes_{p-r}f\|^2_{\HH^{\otimes
2r}}+\|g\otimes_{q-r}g\|^2_{\HH^{\otimes 2r}}\big).\notag
\end{eqnarray}

2. One has the identity:
\begin{eqnarray}
E[F^4]-3E[F^2]^2
&=&
\sum_{r=1}^{p-1} p!^2\binom{p}{r}^2\left\{ 
\|f\otimes_{r} f\|^2_{\HH^{\otimes 2p-2r}}
+\binom{2p-2r}{p-r}\|f\widetilde{\otimes}_{r} f\|^2_{\HH^{\otimes 2p-2r}}
\right\}.\notag\\
\label{gio}
\end{eqnarray}
\end{lemma}
{\it Proof}. 
Inequalities (\ref{ineq1})-(\ref{ineq2}) correspond to \cite[Lemma 3.7]{NPReveillac}
(see also \cite[Proof of Lemma 6]{NuOL}), whereas
identity (\ref{gio}) is shown in \cite[page 182]{NP}. However, for convenience of the reader
(and also because the notation used in \cite{NP} is not exactly the same than our), we
provide here a detailed proof of (\ref{ineq1}), (\ref{ineq2}) and (\ref{gio}).

1. Thanks to the multiplication formula (\ref{multiplication}), we can write
\begin{eqnarray*}
\langle DF,DG\rangle_\HH &=&p\,q\left\langle I_{p-1}(f),I_{q-1}(g)\right\rangle_\HH\\
&=&p\,q\sum_{r=0}^{p\wedge q-1} r!\binom{p-1}{r}\binom{q-1}{r} I_{p+q-2-2r}(f\widetilde{\otimes}_{r+1}g)\\
&=&p\,q \sum_{r=1}^{p\wedge q} (r-1)!\binom{p-1}{r-1}\binom{q-1}{r-1} I_{p+q-2r}(f\widetilde{\otimes}_r g).
\end{eqnarray*}
It follows that
\begin{eqnarray}
&&E\left[\left(\alpha-\frac1q\left\langle DF,DG\right\rangle_\HH\right)^2\right] \label{Murray}\\
&=&\left\lbrace
\begin{array}{l}
\alpha^2+p^2\sum_{r=1}^{p}(r-1)!^2
\binom{p-1}{r-1}^2\binom{q-1}{r-1}^2 (p+q-2r)!
\|f\widetilde{\otimes}_r g\|^2_{\HH^{\otimes (p+q-2r)}} \textrm{ if } p< q,\\\\
\big(\alpha-E[FG]\big)^2+p^2\sum_{r=1}^{p-1}(r-1)!^2 \binom{p-1}{r-1}^4 (2p-2r)!
\|f\widetilde{\otimes}_r g\|^2_{\HH^{\otimes (2p-2r)}} \textrm{ if
} p=q.
\end{array}\notag
\right.
\end{eqnarray}
If $r<p\leq q$ then
\begin{eqnarray}
\|f\widetilde{\otimes}_r g\|^2_{\HH^{\otimes (p+q-2r)}}
&\leq& \|f\otimes_r g\|^2_{\HH^{\otimes (p+q-2r)}}
=\langle f\otimes_{p-r} f, g\otimes_{q-r}g\rangle_{\HH^{\otimes 2r}}\notag\\
&\leq&
\|f\otimes_{p-r}f\|_{\HH^{\otimes 2r}}\|g\otimes_{q-r}g\|_{\HH^{\otimes 2r}}\notag\\
&\leq&\frac12\left(
\|f\otimes_{p-r}f\|_{\HH^{\otimes 2r}}^2+\|g\otimes_{q-r}g\|_{\HH^{\otimes 2r}}^2
\right).\label{aqw1}
\end{eqnarray}
If $r=p<q$, then
\begin{equation}\label{aqw2}
\|f\widetilde{\otimes}_p\, g\|^2_{\HH^{\otimes (q-p)}} \leq
\|f\otimes_p \,g\|^2_{\HH^{\otimes (q-p)}} \leq
\|f\|^2_{\HH^{\otimes p}}\|g\otimes_{q-p}g\|_{\HH^{\otimes 2p}}.
\end{equation}
By plugging these two inequalities into (\ref{Murray}), we deduce both
(\ref{ineq1}) and (\ref{ineq2}).

2. Without loss of generality, in the proof of (\ref{gio}) we can assume that $\HH$ is 
a $L^2$-space of the form $\HH=L^2(A,\mathcal{A},\mu)$.
Let $\sigma$ be a permutation of $\{1,\ldots,2p\}$
(that is, $\sigma\in\mathfrak{S}_{2p}$), and let $f\in\HH^{\odot 2p}$.
If $r\in\{0,\ldots,p\}$ denotes the cardinality of $\{\sigma(1),\ldots,\sigma(p)\}\cap\{1,\ldots,p\}$
then it is readily checked that $r$ is also the cardinality of
$\{\sigma(p+1),\ldots,\sigma(2p)\}\cap\{p+1,\ldots,2p\}$ and that
\begin{eqnarray}
&&\int_{A^{2p}}f(t_1,\ldots,t_p)f(t_{\sigma(1)},\ldots,t_{\sigma(p)})f(t_{p+1},\ldots,t_{2p})
f(t_{\sigma(p+1)},\ldots,t_{\sigma(2p)})d\mu(t_1)\ldots d\mu(t_{2p})\notag\\
&=&\int_{A^{2p-2r}}f\otimes_r f(x_1,\ldots,x_{2p-2r})^2d\mu(x_1)\ldots d\mu(x_{2p-2r}) =
\|f\otimes _r f\|^2_{\HH^{\otimes(2p-2r)}}.\label{ctr}
\end{eqnarray}
Moreover, for any fixed $r\in\{0,\ldots,p\}$, there are $\binom{p}{r}^2(p!)^2$
permutations $\sigma\in\mathfrak{S}_{2p}$ such that
$\#\{\sigma(1),\ldots,\sigma(p)\}\cap\{1,\ldots,p\}=r$.
(Indeed, such a permutation is completely determined by the choice of: $(a)$ $r$ distinct
elements $x_1,\ldots,x_r$ of $\{1,\ldots,p\}$; $(b)$ $p-r$ distinct elements $x_{r+1},\ldots,x_p$
of $\{p+1,\ldots,2p\}$; $(c)$ a bijection between $\{1,\ldots,p\}$ and $\{x_1,\ldots,x_p\}$;
$(d)$ a bijection betwenn $\{p+1,\ldots,2p\}$ and $\{1,\ldots,2p\}\setminus \{x_1,\ldots,x_p\}$.)
Now, observe that the symmetrization of $f\otimes f$ is given by
\[
f\widetilde{\otimes} f(t_1,\ldots,t_{2p}) = \frac{1}{(2p)!}
\sum_{\sigma\in\mathfrak{S}_{2p}} f(t_{\sigma(1)},\ldots,t_{\sigma(p)})
f(t_{\sigma(p+1)},\ldots,t_{\sigma(2p)}).
\]
Therefore,
\begin{eqnarray*}
\|f\widetilde{\otimes} f\|^2_{\HH^{\otimes 2p}}&=&
\frac{1}{(2p)!^2}
\sum_{\sigma,\sigma'\in\mathfrak{S}_{2p}}
\int_{A^{2p}} 
f(t_{\sigma(1)},\ldots,t_{\sigma(p)})f(t_{\sigma(p+1)},\ldots,t_{\sigma(2p)})\\
&&\hskip3cm\times 
f(t_{\sigma'(1)},\ldots,t_{\sigma'(p)})f(t_{\sigma'(p+1)},\ldots,t_{\sigma'(2p)})
d\mu(t_1)\ldots d\mu(t_{2p})\\
&=&
\frac{1}{(2p)!}
\sum_{\sigma\in\mathfrak{S}_{2p}}
\int_{A^{2p}} 
f(t_{1},\ldots,t_{p})f(t_{p+1},\ldots,t_{2p})\\
&&\hskip3cm\times 
f(t_{\sigma(1)},\ldots,t_{\sigma(p)})f(t_{\sigma(p+1)},\ldots,t_{\sigma(2p)})
d\mu(t_1)\ldots d\mu(t_{2p})\\
&=&\frac{1}{(2p)!}\sum_{r=0}^p
\sum_{\substack{\sigma\in\mathfrak{S}_{2p}\\
\{\sigma(1),\ldots,\sigma(p)\}\cap\{1,\ldots,p\}=r
}}
\int_{A^{2p}} 
f(t_{1},\ldots,t_{p})f(t_{p+1},\ldots,t_{2p})\\
&&\hskip3cm\times 
f(t_{\sigma(1)},\ldots,t_{\sigma(p)})f(t_{\sigma(p+1)},\ldots,t_{\sigma(2p)})
d\mu(t_1)\ldots d\mu(t_{2p}).
\end{eqnarray*}
Using (\ref{ctr}), we deduce that
\begin{equation}
(2p)!\|f\widetilde{\otimes} f\|^2_{\HH^{\otimes 2p}}
=2(p!)^2\|f\|_{\HH^{\otimes p}}^4+(p!)^2\sum_{r=1}^{p-1}
\binom{p}{r}^2\|f\otimes_r f\|^2_{\HH^{\otimes(2p-2r)}}.\label{beautyformula}
\end{equation}
On the other hand, we infer from the product formula (\ref{multiplication}) that
\[
F^2=I_p(f)^2=\sum_{r=0}^{p}r!\binom{p}{r}^2 I_{2p-2r}(f\widetilde{\otimes}_r f).
\]
Using the orthogonality and isometry properties of the integrals $I_p$, 
this
yields
\begin{eqnarray*}
E[F^4] &=& \sum_{r=0}^{p} (r!)^2\binom{p}{r}^4 (2p-2r)!
\|f\widetilde{\otimes}_r f\|^2_{\HH^{\otimes (2p-2r)}}\\
&=&(2p)! \|f\widetilde{\otimes} f\|^2_{\HH^{\otimes (2p)}}
+(p!)^2\|f\|^4_{\HH^{\otimes p}}
+\sum_{r=1}^{p-1} (r!)^2\binom{p}{r}^4 (2p-2r)!
\|f\widetilde{\otimes}_r f\|^2_{\HH^{\otimes (2p-2r)}}.
\end{eqnarray*}
By inserting (\ref{beautyformula}) in the previous identity (and because
$(p!)^2\|f\|^4_{\HH^{\otimes p}}=E[F^2]^2$), we
get (\ref{gio}).

\fin

We are now ready to prove Theorem \ref{main1}.
If $Z\in L^4(\Omega)$, as usual we write $\chi_4(Z)=E[Z^4]-3E[Z^2]^2$ for
the fourth cumulant of $Z$.
We deduce from (\ref{gio}) that, for all $p\geq 1$, $f\in\HH^{\odot p}$ and $r\in\{1,\ldots,p-1\}$,
one has $\chi_4(I_p(f))\geq 0$ and 
\[
\|f\otimes_r f\|^2_{\HH^{\otimes 2p-2r}}\leq \frac{r!^2(p-r)!^2}{p!^4}
\,
\chi_4(I_p(f)).
\]
Therefore, if $f,g\in\HH^{\odot p}$, inequality (\ref{ineq1})
leads to
\begin{eqnarray}
E\left[\left(E[I_p(f)I_p(g)]-\frac1p\left\langle DI_p(f),DI_p(g)\right\rangle_\HH\right)^2\right]
&\leq&
\big[ 
\chi_4(I_p(f))
+
\chi_4(I_p(g))
\big]\sum_{r=1}^{p-1}
\frac{r^2(2p-2r)!}{2p^2(p-r)!^2}
\notag\\
&\leq&\frac12	
\big[ 
\chi_4(I_p(f))
+
\chi_4(I_p(g))
\big]\sum_{r=1}^{p-1}
\binom{2r}{r}.\notag\\
\label{azerty1}
\end{eqnarray}
On the other hand, if $p<q$, $f\in\HH^{\odot p}$ and $g\in\HH^{\odot q}$, inequality (\ref{ineq2})
leads to
\begin{eqnarray*}
E\left[\left(\frac1p\left\langle DI_p(f),DI_q(g)\right\rangle_\HH
\right)^2\right]&=&\frac{q^2}{p^2}
E\left[\left(\frac1q\left\langle DI_p(f),DI_q(g)\right\rangle_\HH
\right)^2\right]\\
&\leq& 
E[I_p(f)^2]\sqrt{\chi_4(I_q(g))}+
\frac{1}{2p^2}\sum_{r=1}^{p-1} r^2(p+q-2r)!\\
&&\times\left[
\frac{q!^2}{(q-r)!^2p!^2}\,\chi_4(I_p(f)) +
\frac{p!^2}{(p-r)!^2q!^2}\,\chi_4(I_q(g))
\right]\\
&\leq& 
E[I_p(f)^2]\sqrt{\chi_4(I_q(g))}+
\frac{1}{2}\sum_{r=1}^{p-1} (p+q-2r)!\\
&&\hskip2cm\times\left[
\binom{q}{r}^2\chi_4(I_p(f)) +
\binom{p}{r}^2\chi_4(I_q(g))
\right],
\end{eqnarray*}
so that, if $p\neq q$, $f\in\HH^{\odot p}$ and $g\in\HH^{\odot q}$, one has
that both $E\left[\left(\frac1p\left\langle DI_p(f),DI_q(g)\right\rangle_\HH
\right)^2\right]$ and $E\left[\left(\frac1q\left\langle DI_p(f),DI_q(g)\right\rangle_\HH
\right)^2\right]$ are less or equal than
\begin{eqnarray}
&&E[I_p(f)^2]\sqrt{\chi_4(I_q(g))}+E[I_q(g)^2]
\sqrt{\chi_4(I_p(f))}\label{azerty2}\\
&&\hskip4cm+
\frac{1}2\sum_{r=1}^{p\wedge q-1} (p+q-2r)!\left[
\binom{q}{r}^2\chi_4(I_p(f)) +
\binom{p}{r}^2\chi_4(I_q(g))
\right].\notag
\end{eqnarray}
Since two multiple integrals of different orders are orthogonal, on has that
\[
C_{ij}=E[F_iF_j]=E[I_{q_i}(f_i)I_{q_j}(f_j)]=0\quad\mbox{whenever $q_i\neq q_j$.}\] 
Thus, by using (\ref{azerty1})-(\ref{azerty2}) together with 
$\sqrt{x_1+\ldots+x_n}\leq \sqrt{x_1}+\ldots+\sqrt{x_n}$, we eventually get the desired conclusion
(\ref{3.7}).
\fin

\section{Cumulants for random vectors on the Wiener space}\label{multivariate}
In all this part of the paper, we let the notation of section \ref{S : Malliavin} prevail.
In particular, $X=\{X(h),\,h\in\HH\}$ denotes a given isonormal Gaussian process.

\subsection{Abstract statement}
In this section, by means of the basic operators $D$ and $L$, we calculate
the cumulants of any 
vector-valued functional $F$ of a given isonormal Gaussian process $X$.

First, let us recall the standard multi-index notation. 
A multi-index is a vector $m=(m_1,\ldots,m_d)$
of $\N^d$.
We write
\begin{eqnarray*}
|m|=\sum_{i=1}^d m_i,\quad 
\partial_i =\frac{\partial}{\partial t_i},
\quad \partial^m=\partial_1^{m_1}\ldots\partial_d^{m_d},
\quad x^m = \prod_{i=1}^d x_i^{m_i}.
\end{eqnarray*}
By convention, we have $0^0=1$.
Also, note that $|x^m|=y^m$, where $y_i=|x_i|$ for all $i$.
If $s\in\N^d$, we say that $s\leq m$ if and only if $s_i\leq m_i$ for all $i$.
For any $i=1,\ldots,d$, we let $e_i\in\N^d$ be
the multi-index defined by $(e_i)_j=\delta_{ij}$, with $\delta_{ij}$ the Kronecker symbol.
\begin{defi}
Let $F=(F_1,\ldots,F_d)$ be a $\R^d$-valued random vector such that 
$E|F|^m<\infty$ for some $m\in\N^d\setminus\{0\}$, and let $\phi_F(t)=E[e^{i\langle t,F\rangle_{\R^d}}]$,
$t\in\R^d$, stand for the characteristic function of $F$. The cumulant of order $m$
of $F$ is (well) defined by
\[
\kappa_m(F)=(-i)^{|m|}\partial^m \log \phi_F(t)|_{t=0}.
\]
\end{defi}
For instance, if $F_i,F_j\in L^2(\Omega)$, then $\kappa_{e_i}(F)=E[F_i]$
and $\kappa_{e_i+e_j}(F)={\rm Cov}[F_i,F_j]$.

Now, we need to (recursively) introduce some further notation:
\begin{defi}
Let $F=(F_1,\ldots,F_d)$ be a $\R^d$-valued random vector with
$F_i\in\mathbb{D}^{1,2}$ for each $i$.
Let $l_1,l_2,\ldots$ be a sequence taking values in $\{e_1,\ldots,e_d\}$.
We set $\Gamma_{l_{1}}(F)=F^{l_1}$.
If the random variable $\Gamma_{l_1,\ldots,l_{k}}(F)$ is a
well-defined element of $L^2(\Omega)$ for some $k\geq 1$, we set 
\[
\Gamma_{l_1,\ldots,l_{k+1}}(F)=\langle DF^{l_{k+1}},-DL^{-1}\Gamma_{l_1,\ldots,l_{k}}(F)\rangle_\HH.
\]
\end{defi}
Since
the square-integrability of $\Gamma_{l_1,\ldots,l_k}(F)$ implies that 
$L^{-1}\Gamma_{l_1,\ldots,l_k}(F)\in {\rm Dom}L\subset \mathbb{D}^{1,2}$, 
the definition of $\Gamma_{l_1,\ldots,l_{k+1}}(F)$ makes sense.
 
The next lemma, whose proof is left to the reader because it is an immediate extension
of Lemma 4.2 in \cite{NP-cumulants} to the multivariate case,
gives sufficient conditions on $F$ ensuring that the
random variable $\Gamma_{l_1,\ldots,l_k}(F)$ is a well-defined element of $L^2(\Omega)$.
\begin{lemma}
1. Fix an integer $j\geq 1$, and assume that $F=(F_1,\ldots,F_d)$ is such that $F_i\in\mathbb{D}^{j,2^j}$
for all $i$. 
Let $l_1,l_2,\ldots,l_j$ be a sequence taking values in $\{e_1,\ldots,e_d\}$.
Then, for all $k=1,\ldots,j$, we have that $\Gamma_{l_1,\ldots,l_k}(F)$ is a well-defined element of
$\mathbb{D}^{j-k+1,2^{j-k+1}}$; in particular, one has that $\Gamma_{l_1,\ldots,l_j}(F)\in \mathbb{D}^{1,2}\subset 
L^2(\Omega)$ and that the quantity $E[\Gamma_{l_1,\ldots,l_j}(F)]$ is well-defined and finite.

2. Assume that $F=(F_1,\ldots,F_d)$ is such that $F_i\in\mathbb{D}^\infty$ for all $i$. 
Let $l_1,l_2,\ldots$ be a sequence taking values in $\{e_1,\ldots,e_d\}$.
Then, for all $k\geq 1$, the random variable $\Gamma_{l_1,\ldots,l_k}(F)$ is a well-defined element of
$\mathbb{D}^{\infty}$.
\end{lemma}

We are now ready to state and prove the main result of this section, which is nothing but
the multivariate extension of Theorem 4.3 in \cite{NP-cumulants}.

\begin{thm}\label{representation-cumulants}
Let $m\in\N^d\setminus\{0\}$. Write $m=l_1+\ldots+l_{|m|}$ 
where $l_i\in\{e_1,\ldots,e_d\}$ for each $i$.
(Up to possible permutations of factors, we have existence and uniqueness of this decomposition of $m$.)
Suppose that the random vector $F=(F_1,\ldots,F_d)$ is such that $F_i\in \mathbb{D}^{|m|,2^{|m|}}$
for all $i$.
Then, we have
\begin{equation}\label{rec1}
\kappa_m(F)=(|m|-1)!\,
E\big[\Gamma_{l_1,\ldots,l_{|m|}}(F)\big].
\end{equation}
\end{thm}
\begin{remark}
{\rm
A careful inspection of the forthcoming proof of Theorem \ref{representation-cumulants} shows 
that the
quantity $E\big[\Gamma_{l_1,\ldots,l_{|m|}}(F)\big]$ 
in (\ref{rec1}) 
is actually symmetric with respect to $l_1,\ldots,l_{|m|}$,
that is, 
\[
\forall \sigma\in\mathfrak{S}_{|m|},\quad
E\big[\Gamma_{l_1,\ldots,l_{|m|}}(F)\big]
=E\big[\Gamma_{l_{\sigma(1)},\ldots,l_{\sigma(|m|)}}(F)\big].\] 
}
\end{remark}
{\it Proof of Theorem \ref{representation-cumulants}}. The proof is by induction on $|m|$.
The case $|m|=1$ is clear because $\kappa_{e_j}(F)=E[F_j]=E[\Gamma_{e_j}(F)]$ for all $j$.
Now, assume that (\ref{rec1}) holds for all multi-indices $m\in\N^d$ such that 
$|m|\leq N$, for some $N\geq 1$ fixed, and let us prove that it continues to hold for all the multi-indices $m$ 
verifying $|m|=N+1$. Let $m\in\N^d$ be such that $|m|\leq N$, and fix $j=1,\ldots,d$.
By applying repeatidely (\ref{NPlemma}) and then the chain rule (\ref{chainrule}), we can write
\begin{eqnarray*}
E[F^{m+e_j}]&=&E[F^m\times \Gamma_{e_j}(F)] \\
&=&E[F^m]E[\Gamma_{e_j}(F)]+E[\langle DF^m,-DL^{-1}\Gamma_{e_j}(F)\rangle_\HH]
\\
&=&E[F^m]E[\Gamma_{e_j}(F)]+\sum_{1\leq i_1\leq |m|}E[F^{m-l_{i_1}}\langle DF^{l_{i_1}},-DL^{-1}\Gamma_{e_j}(F)\rangle_\HH]
\\
&=&E[F^m]E[\Gamma_{e_j}(F)]+\sum_{1\leq i_1\leq |m|}E[F^{m-l_{i_1}}\Gamma_{e_j,l_{i_1}}(F)]
\\
&=&E[F^m]E[\Gamma_{e_j}(F)]+\sum_{1\leq i_1\leq |m|}E[F^{m-l_{i_1}}]E[\Gamma_{e_j,l_{i_1}}(F)]
+
\sum_{
\stackrel{1\leq i_1,i_2\leq |m|}{i_1,i_2\,\mbox{\tiny{different}}}
}
E[F^{m-l_{i_1}-l_{i_2}}\Gamma_{e_j,l_{i_1},l_{i_2}}(F)]
\\
&=&\ldots\\
&=&E[F^m]E[\Gamma_{e_j}(F)]+\sum_{1\leq i_1\leq |m|}
E[F^{m-l_{i_1}}]E[\Gamma_{e_j,l_{i_1}}(F)]\\
&&+
\sum_{
\stackrel{1\leq i_1,i_2\leq |m|}{i_1,i_2\,\mbox{\tiny{different}}}
}
E[F^{m-l_{i_1}-l_{i_2}}]E[\Gamma_{e_j,l_{i_1},l_{i_2}}(F)]\\
&&+\ldots+
\sum_{
\stackrel{1\leq i_1,\ldots,i_{|m|-1}\leq |m|}{i_1,\ldots,i_{|m|-1}\,\mbox{\tiny{pairwise different}}}
}
E[F^{m-l_{i_1}-\ldots-l_{i_{|m|-1}}}]E[\Gamma_{e_j,l_{i_1},\ldots,l_{i_{|m|-1}}}(F)] \\
&&+|m|! E[\Gamma_{e_j,l_{i_1},\ldots,l_{i_{|m|}}}(F)]
\end{eqnarray*}
so that, using the induction property,
\begin{eqnarray*}
E[F^{m+e_j}]&=&E[F^m]\frac{1}{0!}\kappa_{e_j}(F)+\sum_{1\leq i_1\leq |m|}
E[F^{m-l_{i_1}}]\frac{1}{1!}\kappa_{e_j+l_{i_1}}(F)\\
&&+
\sum_{
\stackrel{1\leq i_1,i_2\leq |m|}{i_1,i_2\,\mbox{\tiny{different}}}
}
E[F^{m-l_{i_1}-l_{i_2}}]\frac{1}{2!}\kappa_{e_j+l_{i_1}+l_{i_2}}(F)\\
&&+\ldots+
\sum_{
\stackrel{1\leq i_1,\ldots,i_{|m|-1}\leq |m|}{i_1,\ldots,i_{|m|-1}\,\mbox{\tiny{pairwise different}}}
}
E[F^{m-l_{i_1}-\ldots-l_{i_{|m|-1}}}]\frac{1}{(m-1)!}\kappa_{
e_j+l_{i_1}+\ldots+l_{i_{|m|-1}}
}
(F) \\
&&+|m|! E[\Gamma_{e_j,l_{i_1},\ldots,l_{i_{|m|}}}(F)]\\
&=&\sum_{\stackrel{s\leq m}{|s|\leq m-1}}
E[F^{m-s}]\frac{1}{|s|!}\kappa_{e_j+s}(F)\,\#B_s 
+|m|! E[\Gamma_{e_j,l_{i_1},\ldots,l_{i_{|m|}}}(F)].
\end{eqnarray*}
Here,
$B_s$
stands for the set of pairwise different indices $i_1,\ldots,i_{|s|}\in\{1,\ldots,|m|\}$ such that 
$l_{i_1}+\ldots+l_{i_{|s|}}=s$, whereas $\#B_s$ denotes the cardinality of $B_s$. 
Also, let $D_j=\{i=1,\ldots,|m|:\,l_i=e_j\}$ and observe that
$m=(m_1,\ldots,m_d)$ with $m_j=\# D_j$. 
For any $s\leq m$, it is readily checked that
$\# B_s = \binom{m_1}{s_1}\ldots \binom{m_d}{s_d}|s|!$.
(Indeed, to build a multi-index $s=(s_1,\ldots,s_d)$ so that $s\leq m$, one must 
choose $s_1$ indices among the $m_1$ indices of 
$D_1$ up to $s_d$ indices among the $m_d$ indices of $D_d$, and then the order of the factors
in the sum $l_{i_1}+\ldots+l_{i_{|s|}}$.)
Therefore,
\begin{eqnarray*}
E[F^{m+e_j}]&=&\sum_{\stackrel{s\leq m}{|s|\leq m-1}}
\binom{m_1}{s_1}\ldots\binom{m_d}{s_d}
E[F^{m-s}]\,\kappa_{e_j+s}(F)
+|m|! E[\Gamma_{e_j,l_{i_1},\ldots,l_{i_{|m|}}}(F)]\\
&=&\sum_{s\leq m}
\binom{m_1}{s_1}\ldots\binom{m_d}{s_d}
E[F^{m-s}]\,\kappa_{e_j+s}(F)
+|m|! E[\Gamma_{e_j,l_{i_1},\ldots,l_{i_{|m|}}}(F)]-\kappa_{e_j+m}(F)\\
&=&\sum_{s\leq m}
\binom{m_1}{s_1}\ldots\binom{m_d}{s_d}
(-i)^{|m|-|s|}\partial^{m-s}\phi_F(0)\times (-i)^{|s|+1}\partial^{e_j+s}\log\phi_F(0)\\
&&+|m|! E[\Gamma_{e_j,l_{i_1},\ldots,l_{i_{|m|}}}(F)]-\kappa_{e_j+m}(F)\\
&=&(-i)^{|m|+1}\partial^{m}\big(\phi_F\frac{d}{dt_j}\log\phi_F\big)(0)
+|m|! E[\Gamma_{e_j,l_{i_1},\ldots,l_{i_{|m|}}}(F)]-\kappa_{e_j+m}(F)\\
&=&(-i)^{|m|+1}\partial^{m+e_j}\phi_F(0)
+|m|! E[\Gamma_{e_j,l_{i_1},\ldots,l_{i_{|m|}}}(F)]-\kappa_{e_j+m}(F)\\
&=&E[F^{m+e_j}]
+|m|! E[\Gamma_{e_j,l_{i_1},\ldots,l_{i_{|m|}}}(F)]-\kappa_{e_j+m}(F),
\end{eqnarray*}
leading to
\[
|m|! E[\Gamma_{e_j,l_{i_1},\ldots,l_{i_{|m|}}}(F)]=\kappa_{e_j+m},
\]
implying in turn that (\ref{rec1}) holds with $m$ replaced by $m+e_j$.
The proof by induction is concluded.
\fin

\subsection{The case of vector-valued multiple integrals}
We now focus on the calculation of cumulants associated to random vectors
whose component are in a given chaos.
In (\ref{formula-cumulants}) (and in its proof as well), we use the following convention.
For simplicity, we drop the brackets in the writing of 
$f_{\lambda_1}
\widetilde{\otimes}_{r_2}
\ldots
\widetilde{\otimes}_{r_{|m|-1}}
f_{\lambda_{|m|-1}}$, 
by implicitely assuming that
this quantity is defined iteratively from the left to the right. For instance, 
$f\widetilde{\otimes}_{\alpha}g \widetilde{\otimes}_{\beta}h 
\widetilde{\otimes}_{\gamma}k$ actually
means $((f\widetilde{\otimes}_{\alpha}g) \widetilde{\otimes}_{\beta}h) \widetilde{\otimes}_{\gamma}k$.

For convenience, we restate Theorem \ref{thm1.6} (in the more general context of isonormal
Gaussian process).

\begin{thm}\label{thm-pasmal}
Let $m\in\N^d\setminus\{0\}$ with $|m|\geq 3$. Write $m=l_1+\ldots+l_{|m|}$, 
where $l_i\in\{e_1,\ldots,e_d\}$ for each $i$.
(Up to possible permutations of factors, we have existence and uniqueness
of this decomposition of $m$.)
Consider a $\R^d$-valued random vector 
of the form
\begin{eqnarray*}
F=(F_{1},\ldots,F_{d})=\big(I_{q_1}(f_{1}),\ldots,I_{q_d}(f_{d})\big),
\end{eqnarray*}
where each $f_{i}$ belongs to $\HH^{\odot{q_i}}$.
When $l_k=e_j$, we set $\lambda_k=j$, so that
$F^{l_k}=F_{\lambda_k}$ for all $k=1,\ldots,|m|$.
Then:
\begin{equation}\label{formula-cumulants}
\kappa_m(F)= (q_{\lambda_{|m|}})!(|m|-1)!\sum c_{q,l}(r_2,\ldots,r_{|m|-1})
\langle
f_{\lambda_1}\widetilde{\otimes}_{r_2} f_{\lambda_2}\ldots \widetilde{\otimes}_{r_{{|m|-1}}}f_{\lambda_{|m|-1}},
f_{\lambda_{|m|}}
\rangle_{\HH^{\otimes q_{\lambda_{|m|}}}},
\end{equation}
where the sum $\sum$ runs over all collections of integers $r_2,\ldots,r_{|m|-1}$ such that:
\begin{enumerate}
\item[(i)] $1\leq r_i\leq q_{\lambda_i}$ for all $i=2,\ldots,|m|-1$;
\item[(ii)] $r_2+\ldots+r_{|m|-1}=\frac{
q_{\lambda_1}+\ldots+q_{\lambda_{|m|-1}}-q_{\lambda_{|m|}}}{2}$;
\item[(iii)] $r_2 <\frac{q_{\lambda_1}+q_{\lambda_2}}{2}$, $\ldots$, 
$r_2+\ldots +r_{|m|-2}< \frac{q_{\lambda_1}+\ldots+q_{\lambda_{|m|-2}}}2$;
\item[(iv)] $r_3\leq q_{\lambda_1}+q_{\lambda_2}-2r_2$, 
$\ldots$, 
$r_{|m|-1}\leq q_{\lambda_1}+q_{\lambda_{|m|-2}}-2r_2-\ldots-2r_{|m|-2}$;
\end{enumerate}
and where the combinatorial constants $c_{q,l}(r_2,\ldots,r_{s})$ are recursively defined by the
relations
\[
c_{q,l}(r_2)=q_{\lambda_2}(r_2-1)!\binom{q_{\lambda_1}-1}{r_2-1}\binom{q_{\lambda_2}-1}{r_2-1},
\]
and, for $s\geq 3$,
\begin{eqnarray*}
c_{q,l}(r_2,\ldots,r_{s})&=&q_{\lambda_s}(r_{s}-1)!\binom{q_{\lambda_1}+\ldots+q_{\lambda_s}-2r_2-\ldots - 
2r_{s-1}-1}{r_{s}-1}\\
&&\hskip6.5cm\times
\binom{q_{\lambda_s}-1}{r_{s}-1}c_{q,l}(r_2,\ldots,r_{s-1}).
\end{eqnarray*}
\end{thm}
{\it Proof}.
If $f\in\HH^{\odot p}$ and $g\in\HH^{\odot q}$ ($p,q\geq 1$), the multiplication
formula yields
\begin{eqnarray}
\langle DI_p(f),-DL^{-1}I_q(g)\rangle_\HH &=&p\left\langle I_{p-1}(f),I_{q-1}(g)\right\rangle_\HH\notag\\
&=&q\sum_{r=0}^{p\wedge q-1} r!\binom{p-1}{r}\binom{q-1}{r} I_{p+q-2-2r}(f\widetilde{\otimes}_{r+1}g)\notag\\
&=&q \sum_{r=1}^{p\wedge q} (r-1)!\binom{p-1}{r-1}\binom{q-1}{r-1} I_{p+q-2r}(f\widetilde{\otimes}_r g).\label{recrec}
\end{eqnarray}
Thanks to (\ref{recrec}), it is straightforward to prove by induction on $|m|$ that
\begin{eqnarray}
&&\Gamma_{l_1,\ldots,l_{|m|}}(F)\\
&=&
\sum_{r_2=1}^{q_{\lambda_1}\!\wedge q_{\lambda_2}} 
\ldots\sum_{r_{{|m|}}=1}^{[q_{\lambda_1}+\ldots+q_{\lambda_{{|m|}-1}}-2r_2-\ldots-2r_{{|m|}-1}]\wedge q_{\lambda_{|m|}}}
c_{q,l}(r_2,\ldots,r_{{|m|}})
{\bf 1}_{\{r_2< \frac{q_{\lambda_1}+q_{\lambda_2}}2\}}
\ldots 
\notag\\
&&
\times 
{\bf 1}_{\{
r_2+\ldots+r_{{|m|}-1}< \frac{q_{\lambda_1}+\ldots+q_{\lambda_{{|m|}-1}}}2
\}}I_{q_{\lambda_1}+\ldots+q_{\lambda_{|m|}}-2r_2-\ldots-2r_{{|m|}}}\big(
f_{\lambda_1}\widetilde{\otimes}_{r_2} f_{\lambda_2}\ldots \widetilde{\otimes}_{r_{{|m|}}}f_{\lambda_{|m|}}
\big).\notag\\
\label{for}
\end{eqnarray}
Now, let us take the expectation on both sides of (\ref{for}). We get
\begin{eqnarray*}
&&\kappa_{m}(F)\\
&=&(|m|-1)!E[\Gamma_{l_1,\ldots,l_{|m|}}(F)]\\
&=&(|m|-1)!
\sum_{r_2=1}^{q_{\lambda_1}\!\wedge q_{\lambda_2}} 
\ldots\sum_{r_{{|m|}}=1}^{[q_{\lambda_1}+\ldots+q_{\lambda_{{|m|}-1}}-2r_2-\ldots-2r_{{|m|}-1}]\wedge q_{\lambda_{|m|}}}
\!\!c_{q,l}(r_2,\ldots,r_{{|m|}})
{\bf 1}_{\{r_2< \frac{q_{\lambda_1}+q_{\lambda_2}}2\}}
\ldots 
\notag\\
&&
\hskip.3cm\times 
{\bf 1}_{\{
r_2+\ldots+r_{{|m|}-1}< \frac{q_{\lambda_1}+\ldots+q_{\lambda_{{|m|}-1}}}2
\}}
{\bf 1}_{\{
r_2+\ldots+r_{{|m|}}= \frac{q_{\lambda_1}+\ldots+q_{\lambda_{{|m|}}}}2
\}}
\times
f_{\lambda_1}\widetilde{\otimes}_{r_2} f_{\lambda_2}\ldots \widetilde{\otimes}_{r_{{|m|}}}f_{\lambda_{|m|}}
.\notag
\end{eqnarray*}
Observe that, if $2r_2+\ldots+2r_{|m|}= q_{\lambda_1}+\ldots+q_{\lambda_{|m|}}$ and 
$
r_{|m|}\leq q_{\lambda_1}+\ldots+q_{\lambda_{|m|-1}}-2r_2- \ldots - 2r_{|m|-1}$,
then 
\[
2r_{|m|}=q_{\lambda_{|m|}}+\big( 
q_{\lambda_1}+\ldots+q_{\lambda_{|m|-1}}-2r_2- \ldots - 2r_{|m|-1}
\big)\geq q_{\lambda_{|m|}}+r_{|m|},
\]
that is, $r_{|m|}\geq q_{\lambda_{|m|}}$, so that $r_{|m|}=q_{\lambda_{|m|}}$.
Therefore,
\begin{eqnarray*}
&&\kappa_{m}(F)\\
&=&(|m|-1)!
\sum_{r_2=1}^{q_{\lambda_1}\!\wedge q_{\lambda_2}} 
\ldots\!\!\sum_{r_{{|m|}}=1}^{[q_{\lambda_1}+\ldots+q_{\lambda_{{|m|}-1}}-2r_2-\ldots-2r_{{|m|}-1}]\wedge q_{\lambda_{|m|}}}
\!\!c_{q,l}(r_2,\ldots,r_{{|m|}})
{\bf 1}_{\{r_2< \frac{q_{\lambda_1}+q_{\lambda_2}}2\}}
\ldots 
\notag\\
&&
\times 
{\bf 1}_{\{
r_2+\ldots+r_{{|m|}-1}< \frac{q_{\lambda_1}+\ldots+q_{\lambda_{{|m|}-1}}}2
\}}
{\bf 1}_{\{
r_2+\ldots+r_{{|m|}}= \frac{q_{\lambda_1}+\ldots+q_{\lambda_{{|m|}}}}2
\}}\\
&&\times
\langle
f_{\lambda_1}\widetilde{\otimes}_{r_2} f_{\lambda_2}
\ldots \widetilde{\otimes}_{r_{{|m|-1}}}f_{\lambda_{|m|-1}},
f_{\lambda_{|m|}}
\rangle_{\HH^{\otimes q_{\lambda_{|m|}}}}
,\notag
\end{eqnarray*}
which is the announced result,
since $c_{q,l}(r_2,\ldots,r_{|m|-1},q_{\lambda_{|m|}})=(q_{\lambda_{|m|}})!c_{q,l}(r_2,\ldots,r_{|m|-1})$.
\fin

\subsection{Yet another proof of Theorem \ref{npt-thm}}\label{sec4.3}
As a corollary of Theorem \ref{thm-pasmal}, we can now perform yet another proof of 
the implication $(ii)\to(i)$ 
(the only one which is difficult) in Theorem \ref{npt-thm}.
So, let the notation and assumptions of this theorem prevail, suppose that $(ii)$ is in order,
and let us prove that $(i)$ holds.
Applying the method of moments/cumulants, we are left to prove that the cumulants of $F_n$ verify,
for all $m\in\mathbb{N}^d$,
\[
\kappa_m(F_n)\to \kappa_m(N)=\left\{
\begin{array}{cl}
0&\mbox{if $|m|\neq 2$}\\
C_{ij}&\mbox{if $m=e_i+e_j$}
\end{array}
\right.\mbox{ as $n\to\infty$}.
\]
Let $m\in\mathbb{N}^d\setminus\{0\}$.
If $m=e_j$ for some $j$ (that is, if and only if $|m|=1$), we have $\kappa_m(F_n)=E[F_{j,n}]=0$.
If $m=e_i+e_j$ for some $i,j$ (that is, if and only if $|m|=2$), we have
$\kappa_m(F_n)=E[F_{i,n}F_{j,n}]\to C_{ij}$ by assumption (\ref{cij}).
If $|m|\geq 3$, we consider the expression (\ref{formula-cumulants}).
Thanks to (\ref{gio}), from $(ii)$ we deduce that $\|f_{i,n}\otimes_r f_{i,n}\|_{L^2([0,T]^{q_i})}
\to 0$ as $n\to\infty$ for all $i$, whereas, thanks to (\ref{cij}), we deduce
that $q_i!
\|f_{i,n}\|^2_{L^2([0,T]^{q_i})}
=E[F_{i,n}^2]\to C_{ii}$ for all $i$, so that
$\sup_{n\geq 1}\|f_{i,n}\|_{L^2([0,T]^{q_i})}<\infty$ for all $i$.
Let $r_2,\ldots,r_{|n|-1}$ be some integers such that
$(i)$--$(iv)$ in Theorem \ref{thm-pasmal} are satisfied.
In particular, $r_2<\frac{q_{\lambda_1}+q_{\lambda_2}}2$.
From (\ref{aqw1})-(\ref{aqw2}), it comes that 
$\|f_{\lambda_1,n}\widetilde{\otimes}_{r_2} f_{\lambda_2,n}
\|_{L^2([0,T]^{q_{\lambda_1}+q_{\lambda_2}-2r_2})}\to 0$ as $n\to\infty$.
Hence, using Cauchy-Schwarz inequality successively through
\[
\|g\widetilde{\otimes}_r h\|_{L^2([0,T]^{p+q-2r})}\leq
\|g\otimes_r h\|_{L^2([0,T]^{p+q-2r})}\leq
\|g\|_{L^2([0,T]^{p})}\|h\|_{L^2([0,T]^{q})}
\]
whenever $g\in L_s^2([0,T]^{p})$, $h\in L^2_s([0,T]^{q})$
and $r=1,\ldots,p\wedge q$, we get
that
\[
\langle
f_{\lambda_1,n}\widetilde{\otimes}_{r_2} f_{\lambda_2,n}
\ldots \widetilde{\otimes}_{r_{{|m|-1}}}f_{\lambda_{|m|-1},n},
f_{\lambda_{|m|};n}
\rangle_{L^2([0,T]^{q_{\lambda_{|m|}}})}\to 0\quad\mbox{as $n\to\infty$}.
\]
Therefore, $\kappa_m(F_n)\to 0$ as $n\to\infty$ by (\ref{formula-cumulants}).
\fin

\end{document}